\documentclass[letterpaper, 10 pt, conference]{ieeeconf} 

\IEEEoverridecommandlockouts                              % This command is only
 \pdfoutput=1
\usepackage{amsmath, amssymb}
\usepackage{amsfonts}
\usepackage{graphicx}
\usepackage{enumerate}
\usepackage{ifthen}
\usepackage{multicol}
\usepackage{float}
\usepackage{fancyhdr}
\usepackage[usenames, dvipsnames]{color}
\graphicspath{{figures/}}
\usepackage{mathtools}
\usepackage{multicol}
\usepackage{xcolor}
\usepackage{thmtools,thm-restate} % to restate theorems

\usepackage{cite}
\usepackage{svg}

 \usepackage{cancel} %To be able to cross out text
 \usepackage[normalem]{ulem}

\newtheorem{theorem}{Theorem}
\newtheorem{lemma}[theorem]{Lemma}

\newtheorem{example}{Example}

\newtheorem{definition}{Definition}

\newboolean{showcomments}
\setboolean{showcomments}{true}

\newcommand{\bmat}[1]{\begin{bmatrix}
#1
\end{bmatrix}}

\newcommand{\todo}[1]{  \ifthenelse{\boolean{showcomments}}
{\textcolor{ForestGreen}{TO DO:  #1}}{}}
\newcommand{\suggest}[1]{\ifthenelse{\boolean{showcomments}}
{\textcolor{Orange}{(Suggestion: #1)}}{}}
\newcommand{\alain}[1]{\ifthenelse{\boolean{showcomments}}
{\textcolor{Blue}{(Alain says: #1)}}{}}
\newcommand{\jonas}[1]{\ifthenelse{\boolean{showcomments}}
{\textcolor{ForestGreen}{(Jonas says: #1)}}{}}
\newcommand{\kristian}[1]{\ifthenelse{\boolean{showcomments}}
{\textcolor{Blue}{(Kristian says: #1)}}{}}
\newcommand{\emma}[1]{\ifthenelse{\boolean{showcomments}}
{\textcolor{VioletRed}{(Emma says: #1)}}{}}
\newcommand{\ifneeded}[1]{\ifthenelse{\boolean{showcomments}}
{\textcolor{Gray}{#1}}{}}

\newboolean{showedit}
\setboolean{showedit}{true}
\newcommand{\edit}[1]{\ifthenelse{\boolean{showedit}}
{\textcolor{Blue}{#1}}{}}
\newcommand{\draft}[1]{\ifthenelse{\boolean{showedit}}
{\textcolor{gray}{#1}}{}}

\newcommand{\ts}{\textsuperscript}
\newcommand{\mc}{\mathcal}

\usepackage{array}
\newcolumntype{L}[1]{>{\raggedright\let\newline\\\arraybackslash\hspace{0pt}}m{#1}}
\newcolumntype{C}[1]{>{\centering\let\newline\\\arraybackslash\hspace{0pt}}m{#1}}
\newcolumntype{R}[1]{>{\raggedleft\let\newline\\\arraybackslash\hspace{0pt}}m{#1}}

\usepackage{physics} %to enable use of \norm

\usepackage[font = small]{caption}
\usepackage{subcaption}

\usepackage{tikz}
 \usetikzlibrary{plotmarks}
\usepackage{tikz}
\usepackage{pgfplots}
\pgfplotsset{compat=newest}
\usetikzlibrary{patterns}
\usetikzlibrary{decorations.text}
\usepgfplotslibrary{fillbetween}

\usepackage{flushend}

\providecommand{\figref}{}
\renewcommand{\figref}[1]{Fig.~\ref{#1}}
\providecommand{\secref}{}
\renewcommand{\secref}[1]{Section~\ref{#1}}

\usepackage{lipsum}
\usepackage{mathtools}

\title{\LARGE \bf Performance bounds for multi-vehicle networks with local integrators

}

\author{{Jonas Hansson and Emma Tegling} 
 \thanks{The authors are with the Department of Automatic Control and the ELLIIT Strategic Research Area at Lund University, Lund, Sweden. Email: \{\tt\small{jonas.hansson, emma.tegling}\}@control.lth.se}
        \thanks{This work was partially funded by Wallenberg AI, Autonomous Systems and Software Program (WASP) funded by the Knut and Alice Wallenberg Foundation and the Swedish Research Council through Grant 2019-00691. }}

\begin{document}
\maketitle
\begin{abstract}
    In this work, we consider the problem of coordinating a collection of $n$\ts{th}-order integrator systems. The coordination is achieved through the novel serial-consensus design, which can be seen as a method for achieving a stable closed-loop while only using local relative measurements. Earlier work has shown that second-order serial consensus can stabilize a collection of double integrators with scalable performance conditions, independent of the number of agents and topology. In this paper, we generalize these performance results to an arbitrary order $n\geq 1$. The derived performance bound depends on the condition number, measured in the vector-induced maximum matrix norm, of a general diagonalizing matrix. We provide an exact characterization of how a minimal condition number can be achieved. Third-order serial consensus is illustrated through a case study of PI-controlled vehicular formation, where the added integrators are used to mitigate the effect of unmeasured load disturbances. The theoretical results are illustrated through examples.  
\end{abstract}

\section{Introduction}

Control of complex systems is a field with a rich literature. A few notable examples are control of the power grid, district heating/cooling, multi-agent coordination, and transportation networks. The common theme is that rich behavior emerges from the many interconnections between the agents. This work concerns the problem of multi-agent coordination, a problem pioneered by \cite{FaxMurray2004, OlfatiSaber2004, Jadbabaie2003}. The consensus protocol has seen many significant theoretical contributions and generalizations. For instance, the works \cite{Bamieh2012Sep, Patterson2014, Siami2014, Pates2017} have considered the sensitivity to noise, often studied through the so-called coherence. There are many more results, such as resilience, scale-fragility, and scalable performance, to name a few. 

The consensus protocol has also been generalized to achieve more complex coordination. In \cite{Ren2005}, the second-order consensus protocol was introduced. This protocol simultaneously coordinates velocity and position for a network of double integrators. In its simplest form, it can be formulated
$$\ddot{x}=-r_\mathrm{vel}L\dot{x}-r_\mathrm{pos}Lx+u_\mathrm{ref},$$
where $r_\mathrm{vel},r_\mathrm{pos}$ are positive constants and $L$ is the graph Laplacian. Since the seminal work \cite{Ren2005}, the protocol has been extensively studied; for example, in terms of coherence \cite{Patterson2014}; stability \cite{Studli2017}; string stability \cite{studli2017StringConcepts, FENG2019StringDefs}. 

Following the idea of second-order consensus, a corresponding high-order version was considered in \cite{Ren2007}:
$$x^{(n)}=-r_{(n-1)}Lx^{(n-1)}-\dots-r_{(1)}L\dot{x}-r_{(0)}Lx+u_\mathrm{ref},$$
where $r_{(1)}=r_\mathrm{vel}$, $r_{(0)}=r_\mathrm{pos}$. Like the second-order protocol, this model has received some attention. An interesting feature of the high-order consensus protocol is that it lacks \emph{scalable stability} \cite{TEGLING2023}; that is, the weights $r_{(k)}$ need to be tuned for the specific graph Laplacian that is used to avoid instability when the number of agents grows. 

A more recent variation of the high-order consensus protocol is \emph{serial consensus} \cite{Hansson2023CDC}. In the Laplace domain, the serial consensus system is
$$\left(\prod_{k=1}^nsI+p_kL\right)X(s)=U_\mathrm{ref}(s).$$
In the second-order case, the corresponding time-domain representation is
$$\ddot{x}=-(p_1+p_2)L\dot{x}-p_2p_1L^2x +u_\mathrm{ref}.$$
The idea is to design a stable closed-loop system independently of the graph Laplacian to mitigate the lack of scalable stability in the conventional high-order consensus~\cite{Studli2017}. However, in high-order consensus, stability is not the only concern, but to be practical, it is often necessary that the closed-loop system satisfies additional {performance} criteria, a problem also considered in \cite{Macellari2017}. Here, we are in particular concerned with \emph{scalable performance} criteria, meaning that the transient should not grow with the size of the formation. This prevents problems such as string instability, a well-known problem in the vehicle formation literature; two surveys are \cite{studli2017StringConcepts, FENG2019StringDefs}. In our work \cite{hansson2024arxivTCNS}, we showed that second-order serial consensus can be used to achieve scalable performance and thus avoid the issue of string instability, regardless of communication topology.

Here, we continue our work on the serial consensus system. An outstanding question was whether the performance result generalizes to higher-order consensus protocols. We provide a positive answer to this question.

Although second-order consensus systems are partially motivated by the vehicle formation problem and multi-agent coordination, they cannot reject unmeasured load disturbances, which unavoidably occur in a real-world setting. Motivated by the classic solution of introducing integral control to counteract steady-state errors, we illustrate how to introduce integral action through a third-order serial consensus design. We show that this system can mitigate load disturbances with scalable performance in terms of maximum transient state deviations due to arbitrary initial conditions.

The article is structured as follows. \secref{sec:preliminaries} will introduce the mathematical notation and models used throughout the article. Then, our main theorem and lemma are presented in \secref{sec:main}. The results are later adapted in \secref{sec:formation} to the setting of PI-controlled vehicle formations to show the scalable performance thereof. The vehicle formation results are further elaborated through numerical simulations in \secref{sec:numerical}. Finally, we conclude the paper with our conclusions and some discussion in \secref{sec:discussion}.

\section{Preliminaries} \label{sec:preliminaries}
\subsection{Definitions and network model}
We will by $\mc G(\mc V,\mc E)$ denote a graph with $N=|\mc V|$ vertices and with edge set $\mc E\subset \mc V\times \mc V$. 
The graph can also be represented by the weighted adjacency matrix $W$, which satisfies $w_{ij}>0 \iff (j, i)\in \mc V$. The graph is called undirected if $W^T=W$. The graph Laplacian $L$ is defined as $L=D-W$ where $D=\mathrm{diag}(W\mathbf{1})$.

A graph contains a directed spanning tree if, for some node $i$, a path exists to all other vertices $j\in \mc V \setminus i$. If the graph contains a directed spanning tree, then the associated graph Laplacian has a unique zero eigenvalue, and the rest of the eigenvalues have strictly positive real parts.

Throughout this paper we will denote the standard vector $\infty$-norm of $x\in \mathbb C^N$ as $\|x\|_\infty=\max_i |x_i|$. We will reuse the notation for the corresponding induced matrix norm of $C\in \mathbb C^{M\times N} $, that is, $\|C\|_\infty=\sup_{\|x\|_\infty=1} \|C x\|_\infty=\max_i \sum_{j=1}^N |C_{ij}|$. We denote the standard Kronecker product between two matrices  $A\otimes B$.

\subsection{Serial Consensus}
% Serial consensus, what to include? 
% 
We begin by introducing serial consensus.
\begin{definition}[$n\ts{th}$-order serial consensus]
    The following closed-loop system, expressed in the Laplace domain,
    \begin{equation}
        \left(\prod_{k=1}^n sI+L_k\right)X(s)=U_\mathrm{ref}(s)
        \label{eq:def_serial_consensus}
    \end{equation}
    is called the $n\ts{th}$-order serial consensus system.
\end{definition}
This can be compared to the conventional $n\ts{th}$-order consensus system which in the case $n=2$ would be
$$\ddot{x}=-L_1\dot{x}(t)-L_0x(t)+u_\mathrm{ref}(t).$$
The corresponding second-order serial consensus system is
$$\ddot{x}(t)=-(L_2+L_1)\dot{x}-L_2L_1x+u_\mathrm{ref}.$$
The two protocols are similar and can be seen as two natural ways of generalizing the regular continuous-time consensus protocol, that is, $\dot x=-Lx$, preserving the limitation of only utilizing relative feedback. We refer to \cite{Hansson2023CDC} for a more elaborate discussion.

One practical property of the serial consensus is that there exists a simple sufficient condition for stability; under the assumption that each of the graph Laplacians $L_k$ each contains a directed spanning tree, the serial consensus is known to eventually achieve coordination \cite{Hansson2023CDC}. When the same graph Laplacians are used, i.e., $L_k=p_kL$, then some additional performance results can be proven, which will be studied later in this paper. In this case, the serial consensus system can be represented in the following state-space form
$$
    \bmat{\dot x\\ \ddot x\\\vdots\\x^{(n)}}
\hspace{-.5ex}\!=\!\hspace{-.5ex}\bmat{ & \hspace{-2ex}I \\
       &  & \hspace{-2ex}\ddots\\
       &  &  &\hspace{-2ex} I\\
       -a_0L^n&\hspace{-1.5ex}- a_1L^{n-1}& \hspace{-2ex}\dots &\hspace{-1.5ex}-a_{n-1}L }
       \hspace{-.9ex}\bmat{ x\\ \dot x\\ \vdots\\x^{(n-1)}}\hspace{-.6ex}+\hspace{-.6ex}
       \bmat{\\ \\ \\u_\mathrm{ref}}\!\!,
$$
where $a_k$ are identified with the weights $p_k$ through the relation $\prod_{k=1}^n(s+p_k)=\sum_{k=0}^{n} a_ks^k$.
Another state space representation that will be of interest is given by
\begin{equation}
    \underbrace{\bmat{ L^{n-1}\dot x\\  L^{n-2}\ddot x\\\vdots\\x^{(n)}}}_{\dot{\xi}}
\hspace{-.65ex}=\hspace{-.65ex}\underbrace{\bmat{ & \hspace{-2ex}L \\
       &  & \hspace{-2ex}\ddots\\
       &  &  &\hspace{-2ex} L\\
       -a_0L&\hspace{-1.5ex}- a_1L& \hspace{-1ex}\dots &\hspace{-2.ex}-a_{n-1}L }}_{\mc A=A\otimes L}
       \hspace{-.5ex}\underbrace{\bmat{ L^{n-1}x\\ L^{n-2}\dot{x}\\ \vdots\\x^{(n-1)}}}_\xi\hspace{-.5ex}+\hspace{-.5ex}
       \bmat{\\ \\ \\u_\mathrm{ref}}\!\!,
       \label{eq:xistates}
\end{equation}
where 
\begin{equation}\label{eq:A}
    A=\bmat{& 1 \\
          & & \ddots\\
            & & & 1\\
            -a_0&-a_1 &\dots &-a_{n-1}}.
\end{equation}
Here, we note that $A$ is in the controllable canonical form associated with the transfer function
$\prod_{k=1}^n(s+p_k)^{-1}.$

\subsection{Scalable performance}

The performance of a system can be evaluated in multiple ways. Inspired by the setting of vehicle formations where velocity following and relative distance keeping is required, we define the measurements $e_\mathrm{pos}=L(x-d)$ and $e_\mathrm{vel}=\dot{x}-v_\mathrm{ref}\mathbf{1}$ with $d\in\mathbb{R}^n$ and $v_\mathrm{ref}\in \mathbb R$ to be the quantities that should be kept small throughout the dynamic phase. This leads to the following definition.
\begin{definition}[Scalable Performance]
    A formation controller defined over a growing family of graphs $\{\mathcal G_N\}$ that ensures 
\begin{equation}
    \sup_{t\geq 0} \left\|\bmat{e_p(t)\\e_v(t)}\right\|_\infty \leq \alpha \left\|\bmat{e_p(0)\\e_v(0)}\right\|_\infty,
    \label{eq:scalable_performance}
\end{equation}
where $\alpha$ is fixed and independent of the number of agents $N$, is said to achieve scalable performance.  
\end{definition}
\vspace{.5ex}
The choice of the norm here is important since the maximum deviations from an equilibrium determine whether collisions will occur or speed limits will be violated.

This definition cannot be directly applied to a more general system since performance is application-dependent. Instead, we consider a more general definition for the high-order consensus protocols.
\begin{definition}[Scalable $\mathcal{C}$-performance]
    A formation controller for the system $\dot{\xi}=-\mc A\xi+\mc B w$ with $y=\mc C \xi$ defined over a growing family of graphs $\{\mathcal G_N\}$ that ensures 
    $$\sup_{t\geq0}  \|y\|_\infty \leq \alpha_y\|y(0)\|_\infty +\alpha_w\sup_{t\geq0}\|w\|_\infty, $$
    where $\alpha_y,\alpha_w$ are fixed and independent of the number of agents $N$, is said to achieve scalable $\mathcal{C}$-performance. 
\end{definition}
We note that this definition is close to scalable input-to-state stability \cite{Besselink2018}. However, here we omit the restriction of having a convergence rate independent of the number of agents. In the case of vehicle formations, and using the measurements $e_\mathrm{pos}$ and $e_\mathrm{vel}$, scalable $\mathcal{C}$-performance implies that no vehicle will deviate far from the desired formation and desired velocity at any point in time. These, we argue, are the main objectives in this setting, while a uniform convergence rate is secondary and potentially not even feasible due to locality and gain constraints.

\section{Main Results}\label{sec:main}

Our first result shows that the serial consensus will have a bounded transient as a response to an arbitrary initial condition and thus achieves scalable $\mathcal{C}$-performance.

\begin{theorem}\label{thrm:nth_serial_performance}
    The $n$\ts{th}-order serial consensus system
    \begin{equation}        
        \left(\prod_{k=1}^nsI+p_kL\right)X(s)=U(s),
    \end{equation}
    with $U(s)=0$, $p_k>0$ for all $k = 1,\ldots, n$, and $p_i\neq p_j$ for all $i\neq j$,
       % has an error vector
       with the states
       $\xi(t)=\left[(L^{n-1}x)^\top\!\!\!,\: (L^{n-2}\dot{x})^\top\!\!\!,\:(L^{n-3}\ddot{x})^\top\!\!\!,\:\cdots\!,\:(x^{(n-1)})^\top \right]^\top\!\!\!$ satisfy
    $$\sup_{t\geq 0} \|\xi(t)\|_\infty \leq \|S\|_\infty\|S^{-1}\|_\infty \|\xi(0)\|_\infty,$$
    where  $S$ is any invertible matrix that diagonlizes $A$ in \eqref{eq:A}.
\end{theorem}
\vspace{1ex}
\begin{proof}
    A state-space representation of the problem, using the states $\xi$, is given by \eqref{eq:xistates}. The solution to the initial value problem is provided by $\xi(t)=e^{A \otimes Lt}\xi(0).$
    Since all $p_k$ are unique, and that $A$ is a controllable canonical form of the transfer function $\prod_{k=1}^n(s+p_k)^{-1}$, it follows that $A$ has $n$ unique negative eigenvalues, being all $-p_k$. Therefore, a diagonalizing and invertible matrix $S$ exists, such that $A=-SPS^{-1}$. The transient can be bounded through the following sequence:
    \begin{align*}
    \|\xi(t)\|_\infty&=\|e^{A\otimes Lt}\xi(0)\|_\infty \\
     &=\|(S\otimes I)e^{-P\otimes Lt}(S^{-1}\otimes I)\xi(0)\|_\infty\\
     &\leq\|S\|_\infty\|S^{-1}\|_\infty\|\xi(0)\|_\infty.
    \end{align*}
    Here, we used submultiplicativity and that $P\otimes L$ is a graph Laplacian, thus satisfying $\|e^{-P\otimes Lt}\|_\infty \leq1$ for all $t\geq 0$.
\end{proof}
This result shows that the serial consensus will have scalable performance independent of the graph structure and the number of agents. Since the bound holds for any diagonalizing matrix $S$, there is potential for improvement. The following general lemma identifies the optimal diagonalizing $S$, which generates the tightest upper bound. 
\begin{lemma}\label{lem:minnorm}
    If $M\in\mathbb R^{n\times n}$ has $n$ unique eigenvalues, then 
    \begin{equation}
        \min_{M=SDS^{-1}} \|S\|_\infty\|S^{-1}\|_\infty=\|S_*K\|_\infty,
    \end{equation}
     where $K$ and $D$ are diagonal matrices, $S_*$ any matrix satisfying $M=S_*DS_*^{-1}$, and $K$ given by   $K_{ii}=\sum_j |S_*^{-1}|_{ij}.$
\end{lemma}
\vspace{1ex}
\begin{proof}
    Since the eigenvalues of $M$ are unique, the eigenvectors are unique up to scaling and permutation. Thus, given one diagonalizing matrix $S_*$, all others can be obtained through $S=S_*KQ$, where $K$ is invertible and diagonal and $Q$ is a permutation matrix. The minimization problem can be simplified through
    \vspace{-1ex}\begin{align*}
        \min_{M=SDS^{-1}}  &\|S\|_\infty\|S^{-1}\|_\infty\\
        =\min_{K, Q}  &\|S_*K Q \|_\infty\|Q^\top K^{-1}S_*^{-1}\|_\infty\\
        =\min_{K}  &\|S_*K\|_\infty\|K^{-1}S_*^{-1}\|_\infty.
    \end{align*}
    The choice of $Q$ does not impact the size. For $K$, suppose $\hat{K}$ solves the minimization problem. First, we note that $K=a\hat{K}$ for any $a\neq 0$ is also optimal, since 
    $$\|S_*a\hat{K}\|_\infty\|(aK)^{-1}S_*^{-1}\|_\infty=\|S_*\hat{K}\|_\infty\|\hat{K}^{-1}S_*^{-1}\|_\infty.$$
    So without loss of generality, we can assume $\|\hat{K}^{-1}S_*^{-1}\|_\infty=1$. By simply applying the definition of the norm, we know that at least one row $i$ satisfy
    $$\left|\hat{K}_{ii}^{-1}\right|\sum_{k=1}^n|S^{-1}_{i,k}|=\|\hat{K}^{-1}S^{-1}\|_\infty=1.$$
    For all other rows $j\neq i$, we instead have
    $$\left|\hat{K}_{jj}^{-1}\right|\sum_{k=1}^n|S^{-1}_{j,k}|\leq 1.\vspace{-.5ex}$$
    In particular, we note that the choice $K_{jj}=\sum_{k=1}^n|S^{-1}_{j,k}|$ preserves the size of $\|K^{-1}S^{-1}\|_\infty=\|\hat{K}^{-1}S^{-1}\|_\infty$ while also satisfying
    $|K_{jj}|\leq |\hat{K}_{jj}|$. At the same time, we have
    \begin{align*}        
    \|S_*K\|_\infty&=\max_i \sum_{j=1}^n|[S_*]_{i,j}||K_{jj}|\\
    &\leq \max_i \sum_{j=1}^n|[S_*]_{i,j}||\hat{K}_{jj}|=\|S_*\hat{K}\|_\infty.
    \end{align*}
    This shows that WLOG the optimal $\hat{K}$ can be assumed to satisfy $\hat{K}_{ii}=\sum_{j=1}^n|[S^{-1}]_{i,j}|$, and since there are no more degrees of freedom, this choice is also a solution to the posed minimization problem.    
\end{proof}
The lemma suggests a general numerical method for finding the
best upper bound in Theorem~\ref{thrm:nth_serial_performance}; identify a matrix $S_*$ together with its inverse that diagonalize $A$ and then apply our lemma. This is possible since effective tools exist for finding eigenvectors and matrix inverses. We also want
to point out that, to the best of our knowledge, this result
is not found in the linear algebra literature. An analogous
result holds in the case of $\|\cdot\|_1.$

\section{Case Study: PI-Control of Vehicle Formation}\label{sec:formation}
Here, we consider a network of double-integrator systems
\begin{equation}
    \ddot{x}_i=u_i(x).
    \label{eq:double_integrators}
\end{equation}
When restricted to only using relative feedback, a natural choice of controller is a second-order consensus system; both conventional and serial consensus can work. However, independent of the choice of second-order consensus system, the closed loop cannot reject load disturbances without a stationary error. This can be seen by the fixed point conditions $L_0x=w\neq 0$ and $L_2L_1x=w\neq 0 \implies L_1x \neq 0$ respectively. For measurable disturbances and with perfectly known dynamics, it is possible to design a feedforward term to reject any disturbance. Alternatively, an integral feedback term can be incorporated to handle unmeasured disturbances. In \cite{TEGLING2023}, it was shown that distributed third-order conventional consensus may lose stability as more agents are added to the network. On the other hand, serial consensus is stable for any order \cite{Hansson2023CDC}. We will, therefore, focus on the third-order serial consensus to handle potential load disturbances in the system. Following the serial-consensus design, the following feedback law can be utilized
\begin{equation}
\begin{aligned}
    u(x)&=-\underbrace{(p_1+p_2+p_3)}_{a_\mathrm{v}} L\dot{x}\\
    -&\underbrace{(p_1p_2+p_1p_3+p_2p_3)}_{a_\mathrm{p}}L^2(x-d)-\underbrace{(p_1p_2p_3)}_{a_{\mathrm{I}}}Lz\\
    \dot z&=L^2(x-d).
\end{aligned}
\label{eq:serial_controller}
\end{equation}
The closed-loop system can then be analyzed in terms of the states $\dot{x}$, $L(x-d)$, and $z$ which has the form
\begin{equation}
\bmat{\dot{z}\\L\dot{(x-d)}\!\!\\ \ddot x}\!\!=\!\!
\bmat{ 0&L&0\\
       0&0&L\\
       \!-a_{\mathrm{I}}L&\!\!\!-{a_\mathrm{p}}L&\!\!\!-{a_\mathrm{v}}L}
   \!\!\bmat{z\\L(x-d)\\ \dot x}\!\!+\!\!\bmat{0\\ 0\\w}\!\!.
\label{eq:3rdOrderSerial}   
\end{equation}
This is of the same form as described in Theorem~\ref{thrm:nth_serial_performance}, so it will also exhibit scalable performance in terms of the transient deviation from any initial condition.

\subsection{Disturbance rejection}
Since only relative feedback is used, it is impossible to guarantee that any load disturbance can be rejected; in fact, only the disturbances in the image of $L$ will be observable. For instance, the response to the constant disturbance $w=\mathbf{1}$ with all states starting at the origin will have the response
$$\dfrac{1^T\dot{x}(t)}{N}=t.$$
But when the disturbance lies in the image, or if we only consider the relative states, the disturbance will be asymptotically rejected as per the following theorem.
\begin{theorem}\label{thrm:3order_dserial}
    The third-order serial consensus system \eqref{eq:3rdOrderSerial} with $w(t)=Lw_0$, where $w_0$ is constant, and where the graph underlying $L$ contains a directed spanning tree,  will achieve $L(x-d)\to 0$ and $L\dot{x}\to 0$ as $t\to \infty$ and satisfy
    $$\sup_{t\geq 0}\|\xi(t)\|_\infty\leq \alpha_\xi(p)\|\xi(0)\|_\infty +\alpha_w(p)\|w_0\|_\infty,$$
    with the bounds quantified by the constants $\alpha_w(p) =\frac{2}{p_1p_2p_3}\|S_1\|_\infty\|S_1^{-1}[1,0,0]^\top\!\|_\infty\!$ and $\alpha_\xi(p)=\|S_2\|_\infty\|S_2^{-1}\|_\infty$ where $S_1,S_2$ are any matrices that diagonalize $A$.
\end{theorem}
\vspace{1ex}
The proof is found in the appendix.
This theorem shows that the PI-controlled vehicle formation will have scalable performance, and load disturbances will not lead to any stationary error. We remark that $\alpha_w$ and $\alpha_\xi$ can be optimized by following the procedure outlined in Lemma~\ref{lem:minnorm}.

\section{Numerical Examples}\label{sec:numerical}
First, we give a simple example of how to apply Lemma~\ref{lem:minnorm} and proceed with two vehicle formation examples.
\subsection{Minimizing norm}
\begin{example}\label{examp:lemmabound}
        Here, we want to find the smallest upper bound to the transient of the second-order serial consensus system $(sI+p_2L)(sI+p_1L)X=U$. The companion matrix~$A$ associated with $(s+p_1)(s+p_2)$ is
        $$A=\bmat{
            0 & 1\\
            -p_1p_2 & -(p_1+p_2)
        }.$$
        Since we know that its eigenvalues are $-p_1,-p_2$ it is easy to find the eigenvectors: $s_1=[1,-p_1]^T$ and $s_2=[1,-p_2]^T$. We can define a diagonalizing matrix $S$ and the corresponding inverse as
        $$S_*=\bmat{1 & 1\\-p_1 & -p_2},\;S_*^{-1}=\frac{1}{p_1-p_2}\bmat{-p_2 & -1\\p_1 & 1}.$$
        Choosing $K=\frac{1}{|p_1-p_2|}\mathrm{diag}(p_2+1,p_1+1)$ ensures that all rows of $K^{-1}S_*^{-1}$ have absolute sums $1$. Finally, we calculate
        \begin{align*}            
        \|S_*K\|_\infty&=\frac{\max(p_2+p_1+2,p_1p_2+p_1+p_2+p_2p_1)}{|p_1-p_2|}\\&=\frac{p_2+p_1 +2\max(1,p_1p_2)}{|p_1-p_2|}.\end{align*}
        This bound is also identical to the upper bound found in \cite{hansson2024arxivTCNS}. For a specific $A$, it is, of course, possible to perform all the previous steps with numerical methods instead.
\end{example}
\subsection{PI control of vehicle formation}
\begin{figure}
    \centering
    \includegraphics[width=1\linewidth]{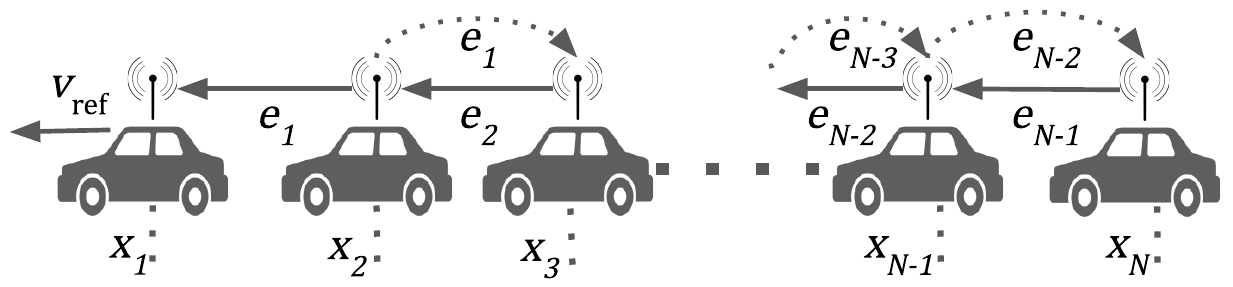}
    \caption{Illustration of vehicle formation control through serial consensus. In this example of a directed string formation, each vehicle measures the relative distance to its neighbor, and message passes its measurements to its follower. }
    \label{fig:cars}
\end{figure}
\begin{example}
    For this example, we will consider a vehicular formation consisting of $N$ identical double integrator systems $\ddot{x}_i=u_i(x)$. A directed path topology where each vehicle only observes its predecessors is used for the control. An illustrating figure for the required communication structure is shown in \figref{fig:cars}. As illustrated, a signaling layer between neighboring vehicles is needed for the agents to use $L^2(x-d)=Le_\mathrm{pos}$ and $Lz$. The associated graph Laplacian is 
    $$L_\mathrm{ahead-path}=\bmat{0 &0\\-1 &1\\& \ddots &\ddots\\
& &-1 &1}.$$
    The serial consensus-based control law \eqref{eq:serial_controller} with $p_1=3$, $p_2=1$, and $p_3=1/3$ and $L=L_\mathrm{ahead-path}$ is used. In \figref{fig:velstep}, a simulation of when the lead vehicle starts at a desired reference velocity while all other vehicles start at rest is shown for $10$ and $40$ vehicles, respectively. A small transient error is observed before the vehicles return to the desired formation as predicted by Theorem~\ref{thrm:nth_serial_performance}.

\begin{figure}
    \centering
    \begin{subfigure}[]{.48\linewidth}
        \centering
        \includegraphics[width=\linewidth]{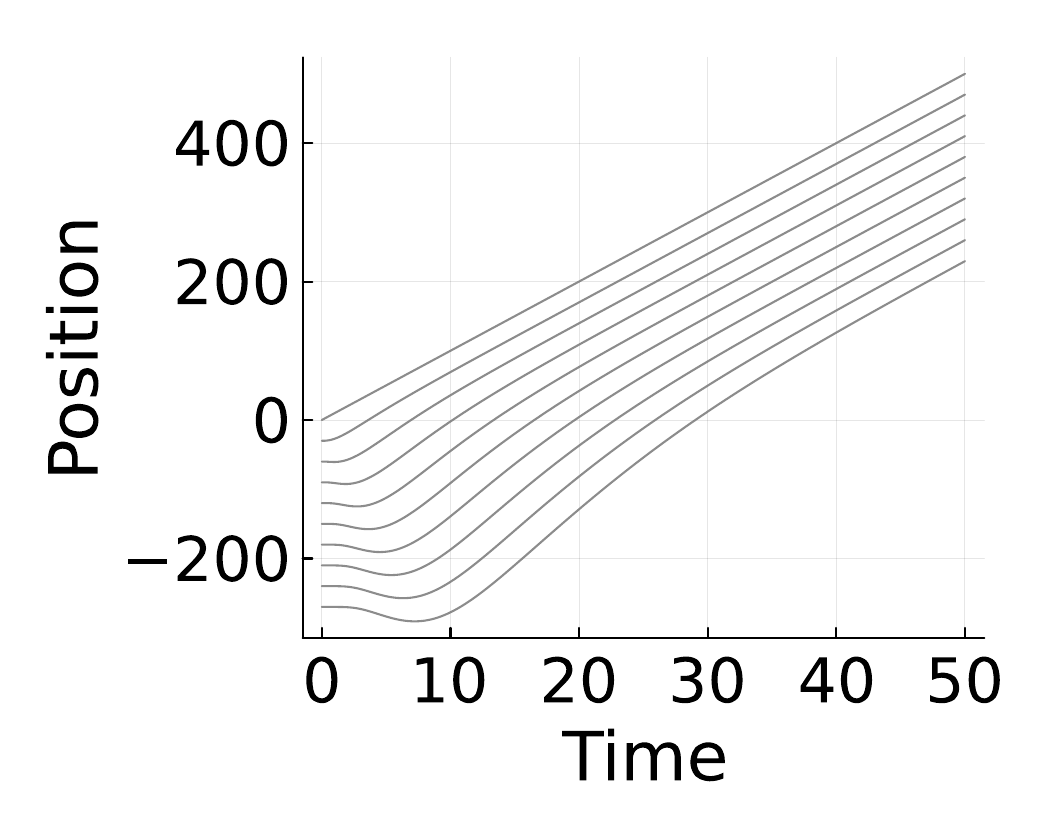}
        \caption{Positions with $N=40$ vehicles.}
        \label{fig:10velstep_pos}
    \end{subfigure} 
    \hfill
        \begin{subfigure}[]{.48\linewidth}
        \centering
        \includegraphics[width=\linewidth]{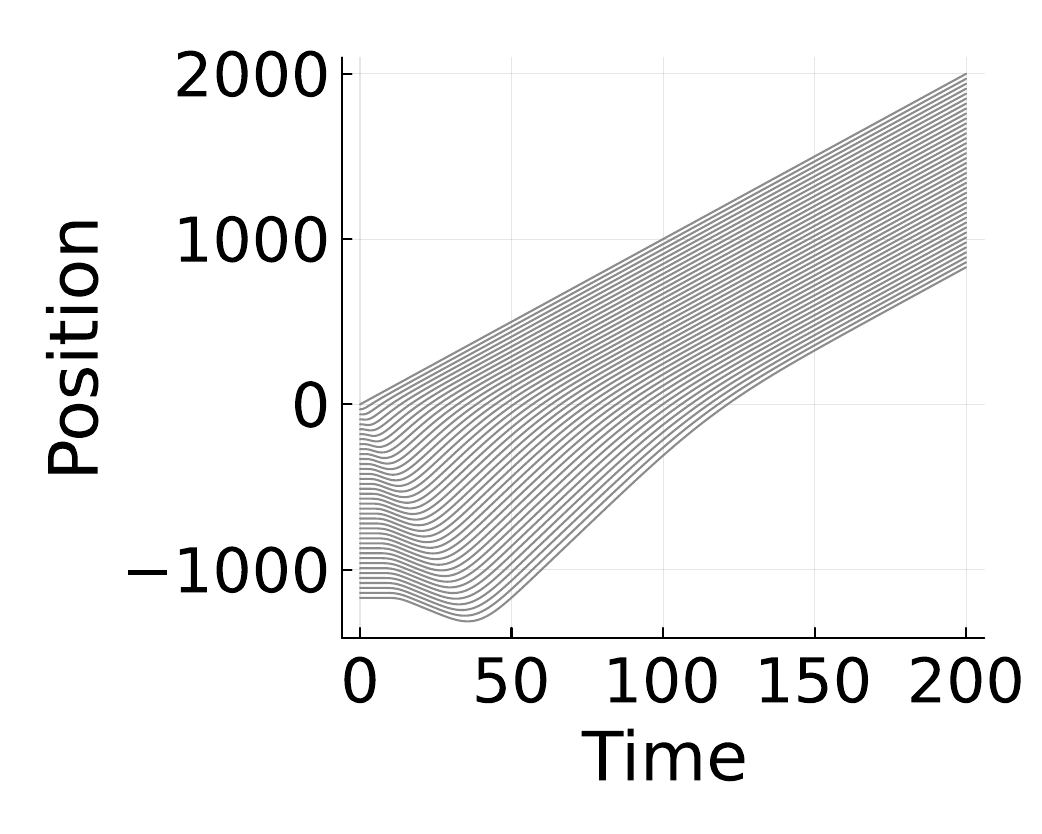}
        \caption{Positions with $N=10$ vehicles.}
        \label{fig:40velstep_pos}
    \end{subfigure} 
   \medskip
    \begin{subfigure}[]{.48\linewidth}
        \centering
        \includegraphics[width=\linewidth]{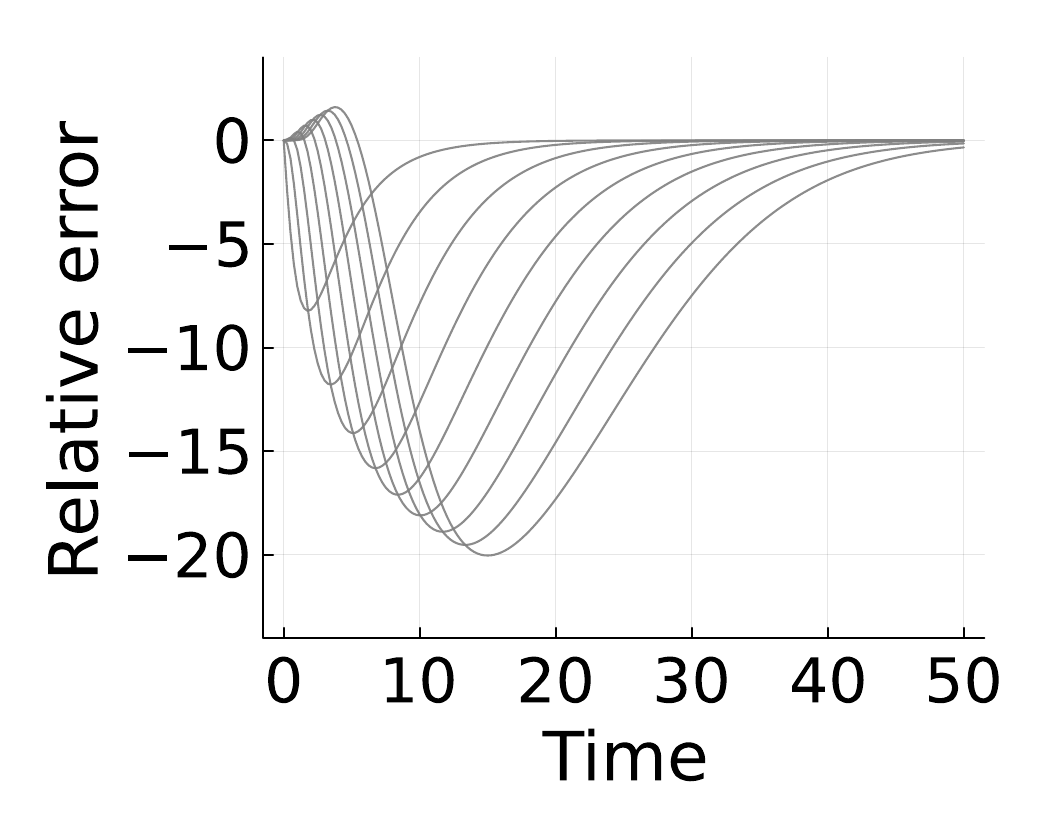}
        \caption{Relative deviations from desired formation with $N=10$. }
        \label{fig:10velstep_relerror}
    \end{subfigure}
    \hfill
    \begin{subfigure}[]{.48\linewidth}
        \centering
        \includegraphics[width=\linewidth]{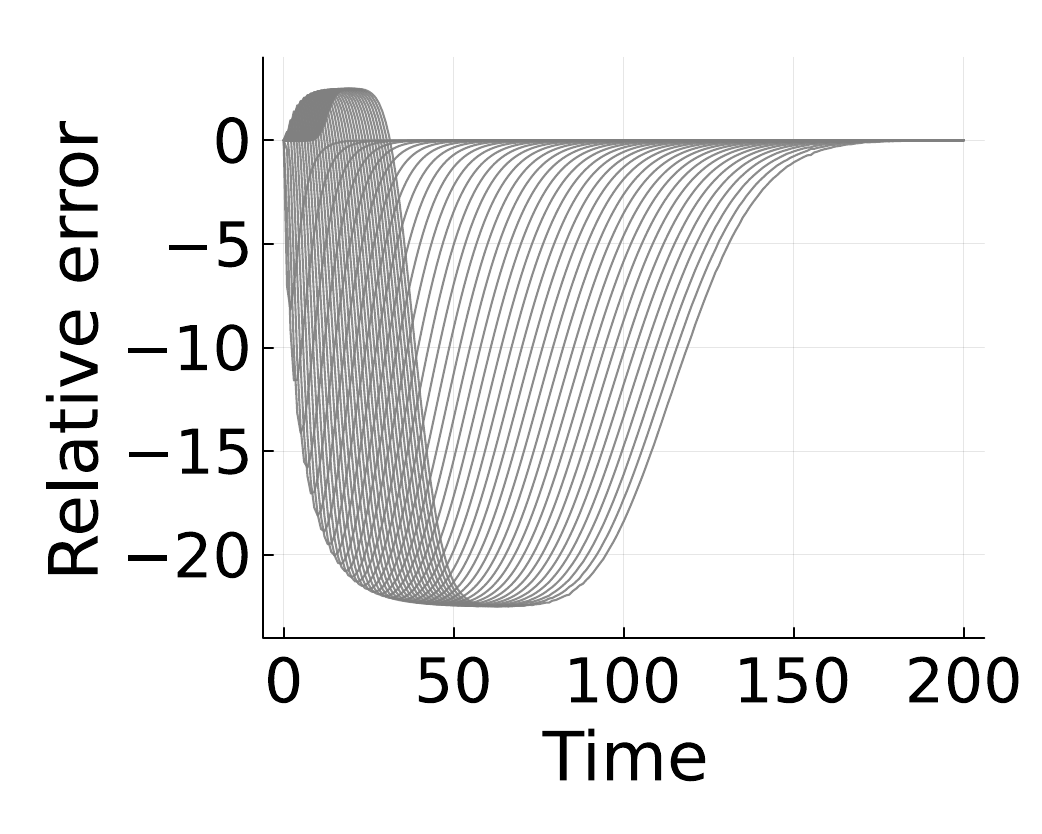}
        \caption{Relative deviations from desired formation with $N=40$. }
        \label{fig:40velstep_relerror}
    \end{subfigure} 
    \caption{Simulation of lead vehicle driving at velocity $v_\mathrm{ref}=10$~m/s and all other vehicles starting from a stand-still. Since the system has scalable performance, the transient error due to any initial condition will be bounded independent of the number of agents. }
    \label{fig:velstep}
\end{figure} 
\end{example}

\subsection{Vehicle formation with load disturbance}
\begin{example}
    In this example, we again consider a vehicle formation consisting of $N$ identical double integrator systems with one virtual leader driving at a constant velocity. The vehicles drive along a road with constant inclination and thus experience a decelerating force (a load disturbance) due to gravity. The individual dynamics are 
    $$\ddot{x}_i=u_i(x,t)-g\dfrac{\theta}{\sqrt{1+\theta^2}},$$
    where $g=9.8$ m/s\ts{2}, is the gravitational acceleration constant and $\theta=0.1$ is the inclination ratio. A second-order serial consensus simulation, with $L=L_\mathrm{ahead-path}$, $p_1=1/3$ and $p_2=3$ is shown in \figref{fig:secondorder_hill_pos} and \figref{fig:secondorder_hill_rel}. Here, a steady-state error exists due to the load disturbance. The third-order case, using the controller \eqref{eq:serial_controller} with $L=L_\mathrm{ahead-path}$, $p_1=1/3$, $p_2=3$, and $p_3=1$ is shown in \figref{fig:3order_hill_pos} and \figref{fig:3order_hill_rel}. In this case, a transient occurs before the vehicles eventually return to the desired formation. Unlike the initial value response, there is no uniform bound independent of the number of agents $N$ for the transient response in this case. In fact, the bound of Theorem~\ref{thrm:3order_dserial} requires that we identify a $w_0$ that solves 
    $$L_\mathrm{ahead-path}w_0\!=\!\dfrac{g\theta}{\sqrt{1+\theta^2}}\bmat{0\\\mathbf{1}}\!\!\implies\!\! w_0\!=\!\mathbf{1}a+\dfrac{g\theta}{\sqrt{1+\theta^2}}\!\bmat{1\\2\\\vdots\\ N}\!\!. $$
    Since this disturbance scales unboundedly with $N$, we cannot give a uniform upper bound for this disturbance.

\begin{figure}
    \centering
    \begin{subfigure}[]{.48\linewidth}
        \centering
        \includegraphics[width=\linewidth]{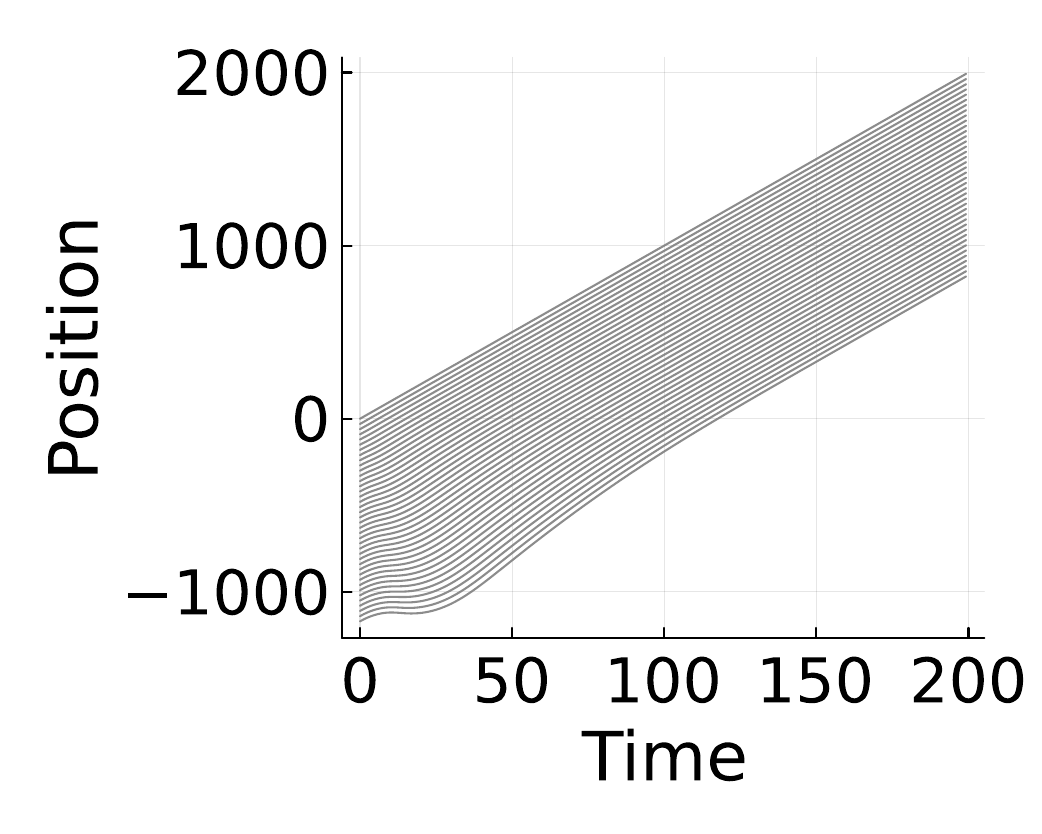}
        \caption{Position for all vehicles, using PI controller. }
        \label{fig:3order_hill_pos}
    \end{subfigure}
    \hfill
    \begin{subfigure}[]{.48\linewidth}
        \centering
        \includegraphics[width=\linewidth]{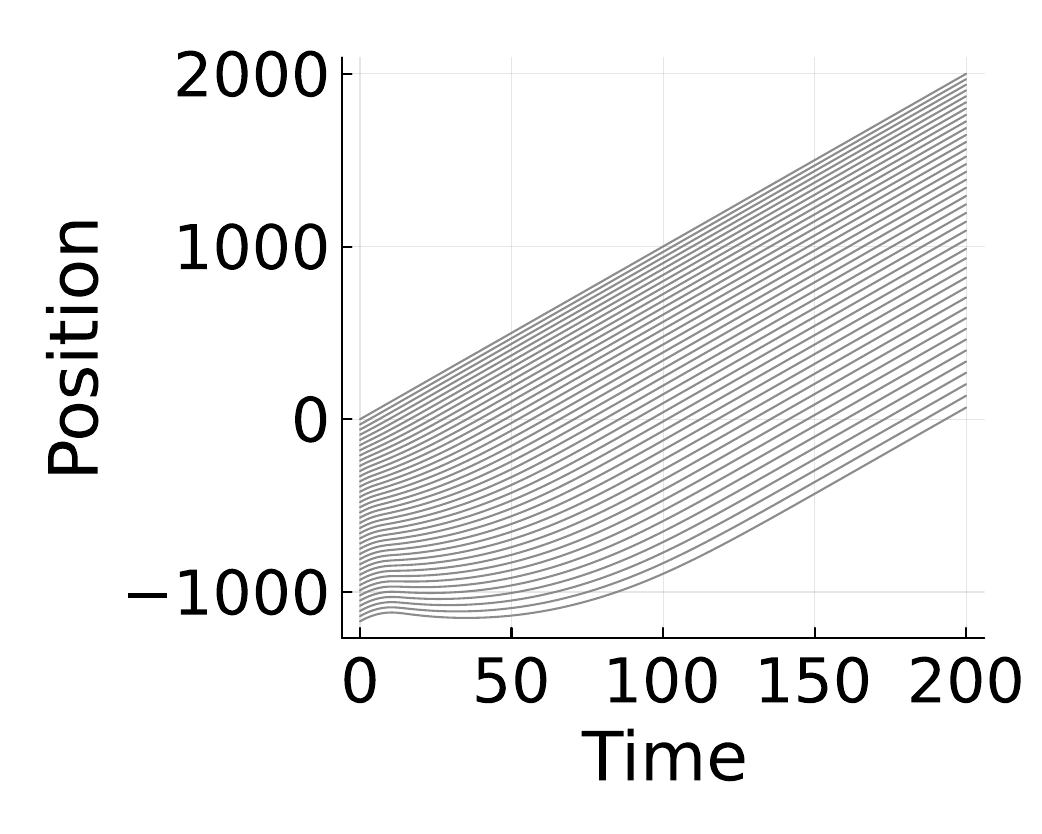}
        \caption{Position for all vehicles, using proportional controller.}
        \label{fig:secondorder_hill_pos}
    \end{subfigure}
   \medskip
    \begin{subfigure}[]{.48\linewidth}
        \centering
        \includegraphics[width=\linewidth]{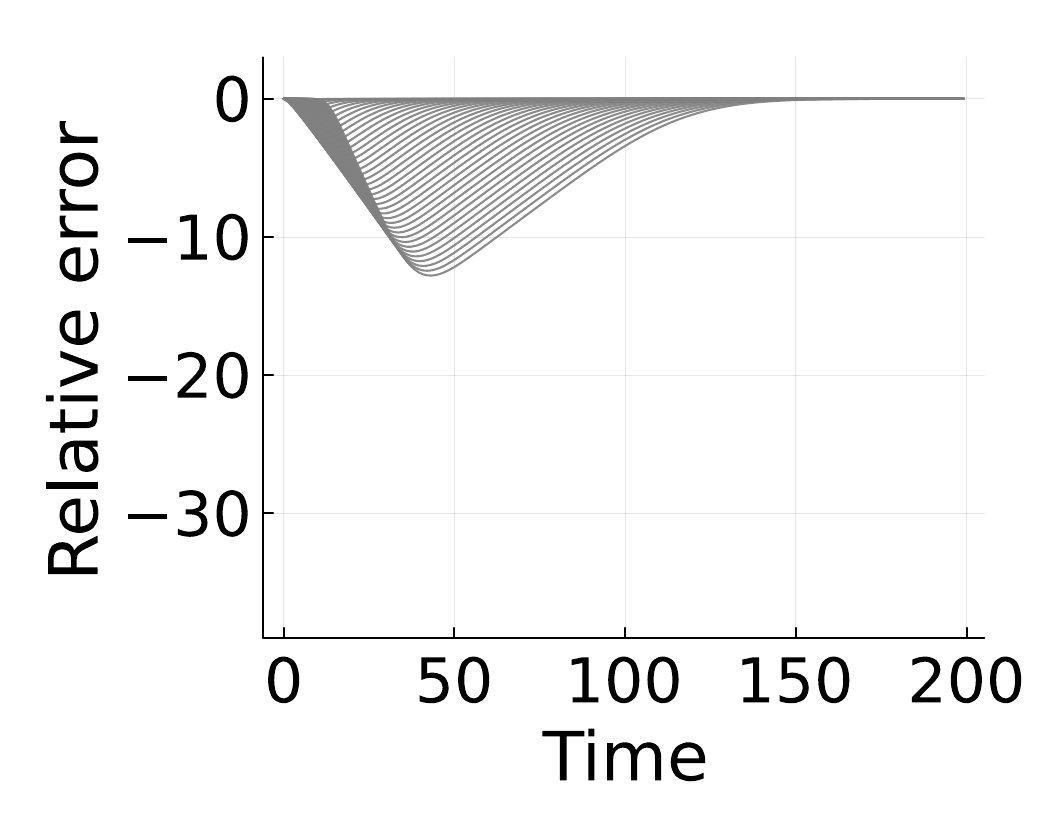}
        \caption{Deviation from desired relative position for all vehicles, using PI controller.}
        \label{fig:3order_hill_rel}
    \end{subfigure} 
    \hfill
        \begin{subfigure}[]{.48\linewidth}
        \centering
        \includegraphics[width=\linewidth]{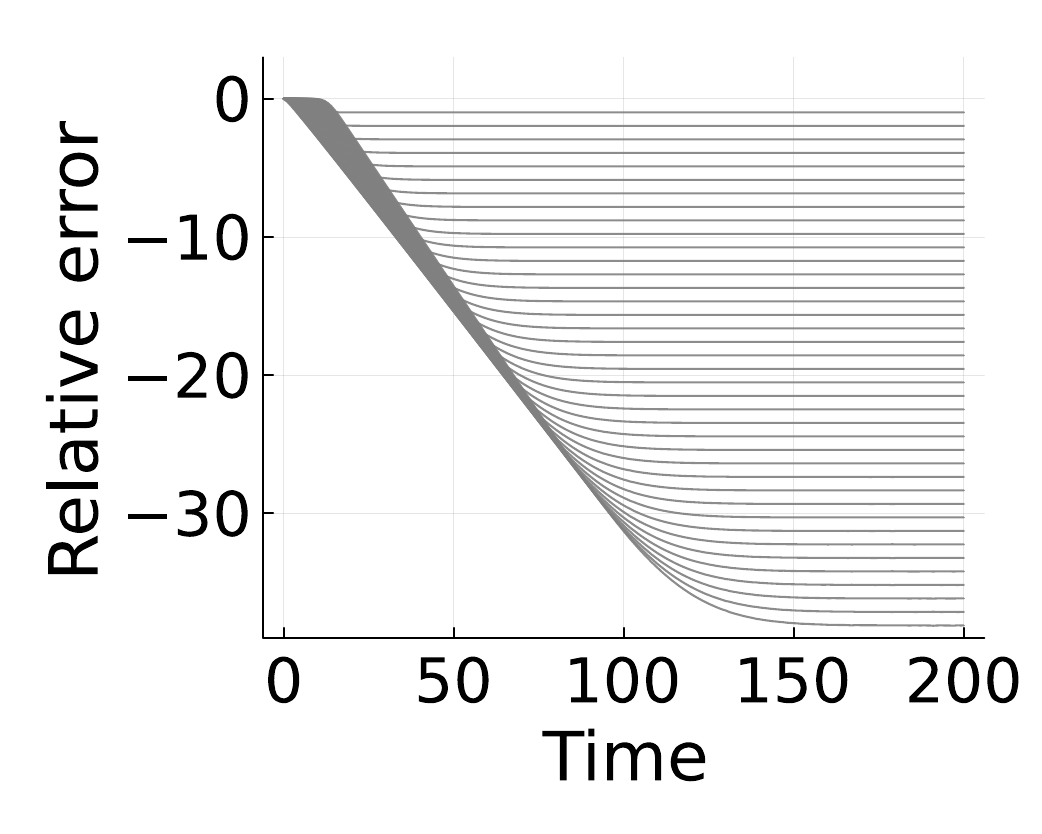}
        \caption{Deviation from desired relative position for all vehicles, using Proportional controller controller.}
        \label{fig:secondorder_hill_rel}
    \end{subfigure} 
    \caption{Simulation of $40$ agents driving uphill at an inclination ratio $\theta=0.1$, and with a leader velocity $10$ m/s. The agents that use a serial consensus-based PI controller experience a transient and then return to the desired spacing. The vehicles controlled with a proportional controller experience a transient before settling with stationary error.}
    \label{fig:3order_hill}
\end{figure}

\end{example}

\section{Discussion and Future Directions}\label{sec:discussion}
This work has expanded the treatment of the newly proposed serial consensus protocol. The derived performance bound also suggests a simple numerical procedure for quantifying its performance.

PI control of vehicle formations serves as an exciting example of how third-order serial consensus can be used to reject load disturbances. Noteworthy, the additional integrator does not fundamentally impact the formation's performance. 

Although it seems that scalable input-to-state stability is not achievable through serial consensus without absolute feedback, a close notion of scalable performance is. We can bound the maximum transient response due to initial conditions and certain load disturbances. 

For future work, it would be interesting to evaluate how measurements such as GPS can be incorporated to achieve stronger notions of scalable performance, such as scalable input-to-state stability. The impact of implementation limitations, such as time delays, should also be addressed.

\section*{Acknowledgement}
We thank Bo Bernhardsson for discussions on Lemma~\ref{lem:minnorm}.

\appendix
\subsection*{Proof of Theorem~\ref{thrm:3order_dserial}}
The solution to \eqref{eq:3rdOrderSerial} is
    $$\xi(t)=e^{A\otimes Lt}\xi(0)+\int_0^te^{A\otimes L(t-\tau)}\left(\bmat{0\\0\\1}\otimes L\right)w_0\,\mathrm{d}\tau,$$
    where $\xi=[z^\top,(L(x-d))^\top,\dot{x}^\top]^\top$. The integral can be simplified through
    \begin{align*}
        \xi_w&=\int_0^te^{A\otimes L(t-\tau)}\left(\bmat{0\\0\\1}\otimes L\right)w_0\,\mathrm{d}\tau\\
        &=\int_0^te^{A\otimes L(t-\tau)}(-A\otimes L)\left(-A^{-1}\bmat{0\\0\\1}\otimes I\right)w_0\,\mathrm{d}\tau\\
        &=(I-e^{A\otimes Lt})\left(-A^{-1}\bmat{0\\0\\1}\otimes I\right)w_0.
    \end{align*}
    The inverse $A^{-1}$ exists and is unique since the eigenvalues are given by the negative constants $-p_1$, $-p_2$, and $-p_3$. It is also easy to show that $-A^{-1}[0,0,1]^T=[1/a_\mathrm{I},0,0]^T$. 
    The complete solution is now 
    $$\bmat{z\\L(x-d)\\ \dot{x}}\!=\!(I-e^{A\otimes Lt})\bmat{w_0/a_\mathrm{I}\\0\\0}+e^{A\otimes Lt}\bmat{z_0\\L(x_0-d)\\\dot{x}_0}\!,
    $$
    Here, we note that the disturbance $w_0$ has a direct term impacting the state $z$ and one part that has the same effect on the states as an initial condition $z_0$ has. Since the consensus equilibrium is stable for any initial condition \cite{Hansson2023CDC}, the effect of $w_0$ will only have a transient effect on $L(x-d)$ and $L\dot{x}$. This proves that $L(x-d)\to0$ and $L(\dot{x})\to 0$
    Clearly, the initial value response bound of Theorem~\ref{thrm:nth_serial_performance} holds for the second term, that is $\alpha_\xi=\|S_1\|_\infty\|S_1^{-1}\|_\infty$, where $S_1$ diagonalizes $A$. The first term can be bounded as follows
    \begin{multline*}
        \|(I-e^{A\otimes Lt})\left([1,0,0]^\top\otimes \frac{w_0}{a_\mathrm{I}}\right)\|_\infty \\
        \leq\|(S_2\otimes I)(I-e^{-P\otimes Lt})(S_2^{-1}[1,0,0]^\top\otimes I)\|_\infty\frac{\|w_0\|_\infty}{a_\mathrm{I}}\\
        \leq \underbrace{\frac{2}{p_1p_2p_3}\|S_2\|_\infty\|S_2^{-1}[1,0,0]^\top\|_\infty}_{\alpha_w}\|w_0\|_\infty,
    \end{multline*}
    where $S_2$ is a potentially different matrix that diagonalizes the matrix~$A$.{\footnotesize $\hfill \blacksquare$}

%\printbibliography 
\bibliographystyle{IEEETran}
\bibliography{references}
\vspace{1mm}
\end{document}